\documentstyle[11pt]{article}
  \newlength{\standardunitlength}
\newtheorem{cor}{Corollary} \newtheorem{lemma}{Lemma}
\newtheorem{theorem}{Theorem} 
\newenvironment{proof}{\noindent {\sc Proof:}}{$\Box$ \vspace{2 ex}}

\begin{document}

\begin{center} {\bf A Probabilistic Approach to Conjugacy Classes in
the Finite Symplectic and Orthogonal Groups} \end{center}

\begin{center}
By Jason Fulman
\end{center}

\begin{center}
Stanford University
\end{center}

\begin{center}
Department of Mathematics
\end{center}

\begin{center}
Building 380, MC 2125
\end{center}

\begin{center}
Stanford, CA 94305, USA
\end{center}

\begin{center}
http://math.stanford.edu/$\sim$ fulman
\end{center}

\begin{center}
fulman@math.stanford.edu
\end{center}

\newpage

\begin{center}
Proposed running head: A Probabilistic Approach
\end{center}

\begin{center}
Please send proofs to:
\end{center}

\begin{center}
Jason Fulman
\end{center}

\begin{center}
Stanford University
\end{center}

\begin{center}
Department of Mathematics
\end{center}

\begin{center}
Building 380, MC 2125
\end{center}

\begin{center}
Stanford, CA 94305, USA
\end{center}

\begin{center}
fulman@math.stanford.edu
\end{center}

\newpage

\begin{abstract} Markov chains are used to give a purely probabilistic
way of understanding the conjugacy classes of the finite symplectic
and orthogonal groups in odd characteristic. As a corollary of these
methods one obtains a probabilistic proof of Steinberg's count of
unipotent matrices and generalizations of formulas of Rudvalis and
Shinoda. \end{abstract}

\begin{center} Key words: Conjugacy class, eigenvalue, random matrix,
classical group, Markov chain.
\end{center}

\section{Introduction} \label{Introduction}

	For compact Lie groups, conjugacy classes are essentially
eigenvalues up to the Weyl group \cite{Ad}. Thus the enormous physics
and mathematical literature on eigenvalues of random matrices (see
\cite{M} for a survey) is often a study of conjugacy classes. Although
the field is rapidly evolving, perhaps the closest thing to a true
probabilistic understanding of eigenvalues of such matrices are the
papers of Dyson \cite{D1,D2}. An equally rich theory exists for the
symmetric groups. The cycles of random permutations have a
probabilistic description using Poisson processes \cite{LS}. The small
cycles of a randomly chosen permutation are asymptotically Poisson,
the medium length cycles relate to Brownian motion \cite{DeP}, and the
long cycles relate to stick breaking \cite{VS}. The cycle structure of
random permutations has numerous applications to real world problems
such as population genetics (for this and more see the reference list
in \cite{F1}).

	Some years back Persi Diaconis observed to the author that a
probabilistic understanding of conjugacy classes of finite groups of
Lie type was missing and urged him to find one. The paper \cite{F2}
provided a useful and beautiful picture for the finite general linear
and unitary groups, with connections to symmetric function theory, but
gave only partial results for the symplectic and orthogonal
groups. The purpose of this article is to complete the program for the
finite classical groups. To this two caveats should be added. First,
only odd characteristic symplectic and orthogonal groups are
considered. As is clear from \cite{W}, the even characteristic
conjugacy classes have a very complicated description; a different
view is given in \cite{FNP}. Second, the current paper lumps together
unipotent conjugacy classes with the same underlying Jordan form. As
indicated at the end of the paper, this can be remedied, but the
resulting formulas seem too complicated to be useful.

	Before describing the contents of this paper, it is worth
remarking that the probabilistic study of Jordan forms of unipotent
upper triangular matrices over a finite field has a fascinating theory
behind it. From the theory of wild quivers there is a provable sense
in which conjugacy classes of upper triangular matrices over a finite
field have no simple description; hence the reduction to Jordan form
is necessary. A lovely probabilistic description of Jordan form is
given in \cite{K} and is exploited in \cite{B}. The survey \cite{F3}
links their work with symmetric function theory and potential theory.

	The main motivation for the current paper is \cite{F4}, which
gave a Markov chain description of the conjugacy classes of the finite
general linear and unitary groups. It is worth recalling the general
nature of that description, as a variation of it occurs here. The
conjugacy classes of $GL(n,q)$ are parameterized by rational canonical
form; for each irreducible polynomial $\phi \neq z$ over $F_q$, one
chooses a partition $\lambda_{\phi}$ of an integer $|\lambda_{\phi}|$,
subject to the constraint that $\sum_{\phi} deg(\phi)
|\lambda_{\phi}|=n$. One can then define a probability measure $M$ on
the set of all partitions of all natural numbers by taking the limit
as $n \rightarrow \infty$ of $\lambda_{z-1}$ for a uniformly chosen
element of $GL(n,q)$. (The polynomial $z-1$ is taken without loss of
generality. For other polynomials one simply replaces $q$ by $q$
raised to the degree of the polynomial in all formulas. Furthermore,
asymptotically the distributions on partitions for different
polynomials are independent). Recall that partitions can be viewed
geometrically. For example the diagram of the partition $(5441)$ is:

\[ \begin{array}{c c c c c}
                . & . & . & . & .  \\
                . & . & . & . &    \\
                . & . & . & . &    \\
                . & & & &  
          \end{array} \] The Markov chain method of sampling from $M$
operates by choosing the size of the first column according to a
certain probability distribution; then given that column $i$ has size
$a$, column $i+1$ will have size $b$ with probability $K(a,b)$. The
remarkable fact is that the probabilities $K(a,b)$ are independent of
$i$, yielding a Markov chain. An immediate consequence of this
viewpoint was an elementary probabilistic proof of the
Rogers-Ramanujan identities, which suggested generalizations to
quivers.

	The main result of this note is that a similar description
occurs for the symplectic and orthogonal cases, except that now the
description will require two Markov chains $K_1$ and $K_2$ defined on
the natural numbers. These Markov chains have the property that they
can never move up, and $K_1$ has the additional property that it can
only move down by an even amount. For the symplectic case, steps with
column number $i$ odd use $K_1$ and steps with column number $i$ even
use $K_2$. For the orthogonal case, steps with column number $i$ odd
use $K_2$ and steps with column number $i$ even use $K_1$. The Markov
chains $K_1,K_2$ are the same for both cases. The only disappointing
aspect of our result is that the product matrices $K_1K_2$ and
$K_2K_1$ do not seem to have a simple diagonalization; this blocked us
from proving Rogers-Ramanujan type identities for the symplectic and
orthogonal groups.

	The structure of this paper is as follows. Section
\ref{recall} recalls the conjugacy classes of the symplectic and
orthogonal groups and defines measures on partitions from them, giving
combinatorially useful rewritings. Section \ref{Markov} begins with
generalizations of formulas of Rudavlis and Shinoda \cite{RS},
\cite{S} and proves the aforementioned description in terms of Markov
chains. In fact it is shown that the set-up extends to a more general
family of measures on partitions with a parameter $u$.

\section{Conjugacy Classes and Measures on Partitions} \label{recall}

	Let $\lambda$ be a partition of some non-negative integer
$|\lambda|$ into parts $\lambda_1 \geq \lambda_2 \geq \cdots$. Let
$m_i(\lambda)$ be the number of parts of $\lambda$ of size $i$, and
let $\lambda'$ be the partition dual to $\lambda$ in the sense that
$\lambda_i' = m_i(\lambda) + m_{i+1}(\lambda) + \cdots$. Let
$n(\lambda)$ be the quantity $\sum_{i \geq 1} (i-1) \lambda_i$. It is
also useful to define the diagram associated to $\lambda$ as the set
of points $(i,j) \in Z^2$ such that $1 \leq j \leq \lambda_i$. We use
the convention that the row index $i$ increases as one goes downward
and the column index $j$ increases as one goes to the right. So the diagram
of the partition $(5441)$ is as in the introduction.

	The following combinatorial lemma about parititions will be
helpful in what follows. For a proof, one simply uses the fact that
$\lambda_i'= m_i(\lambda)+m_{i+1}(\lambda)+\cdots$.

\begin{lemma} \label{combinatorics}
\[ \sum_{h<i} 2hm_h(\lambda)m_i(\lambda) + \sum_i (i-1) m_i(\lambda)^2
= \sum_i (\lambda_i')^2 - \sum_i m_i(\lambda)^2\]
\end{lemma}

	We also recall the following formulas for the sizes of finite
symplectic and orthogonal groups in odd characteristic. 

\begin{eqnarray*}
|Sp(2n,q)| & = & q^{n^2} \prod_{i=1}^n (q^{2i}-1)\\
|O^{\pm}(2n+1,q)| & = & 2q^{n^2} \prod_{i=1}^n (q^{2i}-1)\\
|O^{\pm}(2n,q)| & = & 2q^{n^2-n} (q^n \mp 1) \prod_{i=1}^{n-1} (q^{2i}-1)
\end{eqnarray*} 

\subsection{Symplectic groups} 

	Wall \cite{W} parametrized the conjugacy classes of the finite
symplectic groups and found formulas for their sizes. Let us recall
his parametrization for the case of odd characteristic. Given a a
polynomial $\phi(z)$ with coefficients in $F_q$ and non vanishing
constant term, define a polynomial $\bar{\phi}$ by

\[ \bar{\phi} = \frac{z^{deg(\phi)} \phi(\frac{1}{z})}{\phi(0)}. \]
Wall showed that a conjugacy class of $Sp(2n,q)$ corresponds to the
following data. To each monic, non-constant, irreducible polynomial
$\phi \neq z \pm 1$ associate a partition $\lambda_{\phi}$ of some
non-negative integer $|\lambda_{\phi}|$. To $\phi$ equal to $z-1$ or
$z+1$ associate a symplectic signed partition $\lambda(\pm)_{\phi}$,
by which is meant a partition of some natural number
$|\lambda(\pm)_{\phi}|$ such that the odd parts have even
multiplicity, together with a choice of sign for the set of parts of
size $i$ for each even $i>0$.

\begin{center}
Example of a Symplectic Signed Partition
\end{center}
	
\[ \begin{array}{c c c c c c c c}
	&  . & . & . & . & .  \\
	&  . & . & . & . & .  \\
	+  & . & . & . & . &  \\
	&  . & . & . &  &  \\
	&  . & . & . &  &   \\
	-  & . & . & & &   \\
	 & . & . & & &  
	  \end{array} \] Here the $+$ corresponds to the parts of size
4 and the $-$ corresponds to the parts of size 2. This data represents
a conjugacy class of $Sp(2n,q)$ if and only if:

\begin{enumerate}

\item $|\lambda_{z}|=0$
\item $\lambda_{\phi}=\lambda_{\bar{\phi}}$
\item $\sum_{\phi=z \pm 1} |\lambda(\pm)_{\phi}| + \sum_{\phi \neq z \pm 1} |\lambda_{\phi}| deg(\phi)=2n$
\end{enumerate}

	Let

\[ A_{Sp}(\phi^{i}) =             \left\{ \begin{array}{ll}
																																			 |Sp(m_i(\lambda(\pm)_{\phi}),q)|
& \mbox {if i odd,} \ \phi = z \pm 1\\
																																		 q^{\frac{m_i(\lambda(\pm)_{\phi})}{2}}
|O(m_i(\lambda(\pm)_{\phi}),q)|	 & \mbox{if
i even,} \ \phi = z \pm 1\\

 |U(m_i(\lambda_{\phi}),q^{\frac{deg(\phi)}{2}})| & \mbox{if} \
\phi=\bar{\phi} \neq z \pm 1\\
					
 |GL(m_i(\lambda_{\phi}),q^{deg(\phi)})|^{\frac{1}{2}} &\mbox{if} \ \phi \neq
\bar{\phi}.																															\end{array}
			\right.			 \] where $O(m_i(\lambda_{\phi}),q)$ is the orthogonal group with
the same sign as the sign associated to the parts of size $i$.

	Theorem \ref{distSp} is implicit in the discussion in
\cite{F1}. The three ingredients in its proof are Wall's formulas for
conjugacy class sizes \cite{W}, the deduction that the cycle index of
the symplectic groups factors, and the fact that the formulas in the
statement of Theorem \ref{distSp} define probability measures
(i.e. the asserted probabilities sum to one). This third fact will be
deduced in the proof of Theorem \ref{MarkovSp}, using only an identity
of Cauchy. It is worth emphasizing that neither Steinberg's count of
unipotent elements nor the formulas of Rudvalis and Shinoda in Section
\ref{Markov} are needed to prove the third fact. 

\begin{theorem} \label{distSp} Fix some value of $u$ with
$0<u<1$. Then pick a non-negative even integer with the probability of
getting $2n$ equal to $(1-u^2)u^{2n}$ and pick uniformly in
$Sp(2n,q)$. Let $\Lambda(\pm)_{z-1}, \Lambda(\pm)_{z+1},
\Lambda_{\phi}$ be the random variables corresponding to the conjugacy
class data of the chosen element of $Sp(2n,q)$. Then, aside from the
fact that $\Lambda_{\phi}= \Lambda_{\bar{\phi}}$, any finite number of
these random variables are independent, with probability laws

\begin{eqnarray*}
&&Prob(\Lambda(\pm)_{z-1}=\lambda(\pm)_{z-1})\\
 &=& \frac{ \prod_{r=1}^{\infty}
(1-u^2/q^{2r-1}) u^{|\lambda(\pm)_{z-1}|}}{ q^{[\sum_{h<i}
hm_h(\lambda(\pm)_{z-1})m_i(\lambda(\pm)_{z-1}) + \frac{1}{2} \sum_i
(i-1)m_i(\lambda(\pm)_{z-1})^2]} \prod_i A_{Sp}((z-1)^i)} \\
&&Prob(\Lambda(\pm)_{z+1}=\lambda(\pm)_{z+1})\\
 &=& \frac{ \prod_{r=1}^{\infty}
(1-u^2/q^{2r-1}) u^{|\lambda(\pm)_{z+1}|}}{q^{[\sum_{h<i}
hm_h(\lambda(\pm)_{z+1})m_i(\lambda(\pm)_{z+1}) + \frac{1}{2} \sum_i
(i-1)m_i(\lambda(\pm)_{z+1})^2]} \prod_i A_{Sp}((z+1)^i)} \\
Prob(\Lambda_{\phi}=\lambda_{\phi})
 & = & \frac{ \prod_{r=1}^{\infty}
(1+(-1)^r \frac{u^{deg(\phi)}}{q^{deg(\phi)r/2}}) u^{deg(\phi)
\cdot |\lambda_{\phi}|}}{q^{deg(\phi) [\sum_{h<i}
hm_h(\lambda_{\phi})m_i(\lambda_{\phi}) + \frac{1}{2} \sum_i
(i-1)m_i(\lambda_{\phi})^2]} \prod_i A_{Sp}(\phi^i)} \ if \ 
\phi=\bar{\phi} \neq z \pm 1\\
Prob(\Lambda_{\phi}=\lambda_{\phi})
 & = & \frac{ \prod_{r=1}^{\infty}
(1- \frac{u^{2 deg(\phi)}}{q^{deg(\phi)r}})u^{2 deg(\phi)
\cdot |\lambda_{\phi}|}}{q^{2 deg(\phi) [ \sum_{h<i}
hm_h(\lambda_{\phi})m_i(\lambda_{\phi}) + \frac{1}{2} \sum_i
(i-1)m_i(\lambda_{\phi})^2]} \prod_i A_{Sp}(\phi^i)} \ if \
\phi \neq \bar{\phi}.
\end{eqnarray*} Furthermore, setting $u=1$ in these formulas yields the laws arising
from the $n \rightarrow \infty$ limit of conjugacy classes of a
uniformly chosen element of $Sp(2n,q)$, and the random variables
corresponding to different polynomials are independent, up to the fact
that $\Lambda_{\phi}= \Lambda_{\bar{\phi}}$. \end{theorem}

	From Theorem \ref{distSp}, one sees that if $\phi=\bar{\phi}$,
then the corresponding measures on partitions are specializations of
those for the unitary groups treated in \cite{F2}. Similarly, if $\phi
\neq \bar{\phi}$, then the corresponding measures on partitions are
specializations of those for the general linear groups treated in
\cite{F2}. As the formulas for $z \pm 1$ are the same, for the rest of
this paper only the partition corresponding to $z-1$ will be studied.

	Combining Theorem \ref{distSp} with Lemma \ref{combinatorics}
leads one to the following measure on symplectic signed partitions:

\[ M^{\pm}_{Sp,u}(\lambda(\pm)) = \prod_{r=1}^{\infty} (1-u^2/q^{2r-1})
\frac{u^{|\lambda(\pm)|}}{ q^{1/2 [\sum_i (\lambda(\pm)_i')^2 -
\sum_i m_i(\lambda(\pm))^2 ]} \prod_i A_{Sp}((z-1)^i)}. \] Forgetting
about signs (i.e. lumping together some conjugacy classes) yields a
measure on underlying shapes which will be denoted by
$M_{Sp.u}$. Using the formulas for the sizes of the finite symplectic
and orthogonal groups given at the beginning of this section, one
arrives at the expression:

\begin{eqnarray*}
M_{Sp,u}(\lambda) & = &
\frac{ \prod_{r=1}^{\infty} (1-u^2/q^{2r-1}) u^{|\lambda|}}{q^{1/2 [\sum_i (\lambda_i')^2 - \sum_i m_i^2 ]}
\prod_{i=1 \ mod \ 2} (q^{\frac{m_i^2}{4}} \prod_{l=1}^{m_i/2}
(q^{2l}-1))}\\
&  & \frac{1}{\prod_{i=0 \ mod \ 2 \atop m_i=0 \ mod \ 2}
(q^{\frac{m_i^2}{4}-\frac{m_i}{2}} \prod_{l=1}^{m_i/2} (q^{2l}-1))
\prod_{i=0 \ mod \ 2 \atop m_i=1 \ mod \ 2} (q^{\frac{m_i^2+1}{4}}
\prod_{l=1}^{(m_i-1)/2} (q^{2l}-1))}
\end{eqnarray*}

\subsection{Orthogonal groups} 

	Wall \cite{W} parametrized the conjugacy classes of the finite
orthogonal groups and found formulas for their sizes. Let us recall
his parametrization for the case of odd characteristic. To each monic,
non-constant, irreducible polynomial $\phi \neq z \pm 1$ associate a
partition $\lambda_{\phi}$ of some non-negative integer
$|\lambda_{\phi}|$. To $\phi$ equal to $z-1$ or $z+1$ associate an
orthogonal signed partition $\lambda(\pm)_{\phi}$, by which is meant a
partition of some natural number $|\lambda(\pm)_{\phi}|$ such that all
even parts have even multiplicity, and all odd $i>0$ have a choice of
sign. For $\phi= z-1$ or $\phi = z+1$ and odd $i>0$, we denote by
$\Theta_i (\lambda(\pm)_{\phi})$ the Witt type of the orthogonal group
on a vector space of dimension $m_i(\lambda(\pm)_{\phi})$ and sign the
choice of sign for $i$.

\begin{center}
Example of an Orthogonal Signed Partition
\end{center}
	
\[ \begin{array}{c c c c c}
	&  . & . & . & .   \\
	&  . & . & . & .   \\
	-  & . & . & . &   \\
	&  . & . &  &   \\
	&  . & . &  &    \\
	+  & . &  & &   \\
	 & . &  & &   
	  \end{array}  \]

	Here the $-$ corresponds to the part of size 3 and the $+$
corresponds to the parts of size 1. The data $\lambda(\pm)_{z-1},
\lambda(\pm)_{z+1}, \lambda_{\phi}$ represents a conjugacy class of
some orthogonal group if:

\begin{enumerate}
\item $|\lambda_{z}|=0$
\item $\lambda_{\phi}=\lambda_{\bar{\phi}}$
\item $\sum_{\phi=z \pm 1} |\lambda(\pm)_{\phi}| + \sum_{\phi \neq z \pm
1} |\lambda_{\phi}| deg(\phi)=n$.
\end{enumerate} In this case, the data represents the conjugacy class of
exactly 1 orthogonal group $O(n,q)$, with sign determined by the
condition that the group arises as the stabilizer of a form of Witt
type:

\[ \sum_{\phi=z \pm 1} \sum_{i \ odd} \Theta_i(\lambda(\pm)_{\phi}) +
\sum_{\phi \neq z \pm 1} \sum_{i \geq 1} i m_i(\lambda_{\phi}) {\bf
\omega}, \] where $\omega$ is the Witt type of the quadratic form
$x^2-\delta y^2$ with $\delta$ a fixed non-square in $F_q$.

	Let

\[ A_{O}(\phi^{i}) =             \left\{ \begin{array}{ll}
												q^{-m_i(\lambda(\pm)_{\phi})/2} |Sp(m_i(\lambda(\pm)_{\phi}),q)|																			
& \mbox{if i even,} \ \phi = z \pm 1\\
																																|O(m_i(\lambda(\pm)_{\phi}),q)|	 & \mbox{if
i odd,} \ \phi = z \pm 1\\

 |U(m_i(\lambda_{\phi}),q^{\frac{deg(\phi)}{2}})| & \mbox{if} \ 
\phi=\bar{\phi} \neq z \pm 1\\
					
 |GL(m_i(\lambda_{\phi}),q^{deg(\phi)})|^{\frac{1}{2}} & \mbox{if} \ \phi \neq
\bar{\phi}.																															\end{array}
			\right.			 \] where $O(m_i(\lambda_{\phi}),q)$ is the orthogonal group with
the same sign as the sign associated to the parts of size $i$.

	Theorem \ref{distO} is implicit in \cite{F1}. The three
ingredients in its proof are Wall's formulas for conjugacy class sizes
\cite{W}, the deduction that the cycle index for the sum of $+$ and
$-$ types of the orthogonal groups factors, and the fact that the
formulas in the statement of Theorem \ref{distO} define probability
measures. As in the symplectic case, this third fact will be deduced
in the proof of Theorem \ref{MarkovO}, using only an identity of Cauchy.

\begin{theorem} \label{distO} Fix some value of $u$ with $0<u<1$. Then
pick a non-negative integer with the probability of getting $0$ equal
to $\frac{(1-u)}{(1+u)}$ and probability of getting $n>0$ equal to
$\frac{2u^n(1-u)}{1+u}$. Choose either $O^+(n,q)$ or $O^-(n,q)$ with
probability $\frac{1}{2}$. Finally select an element uniformly within
the chosen orthogonal group and let $\Lambda(\pm)_{z-1},
\Lambda(\pm)_{z+1}, \Lambda_{\phi}$ be the random variables
corresponding to its conjugacy class data. Then aside from the fact
that $\Lambda_{\phi}= \Lambda_{\bar{\phi}}$, any finite number of
these random variables are independent, with probability laws

\begin{eqnarray*}
&&Prob(\Lambda(\pm)_{z-1}=\lambda(\pm)_{z-1})\\
 &=& \frac{ \prod_{r=1}^{\infty}
(1-u^2/q^{2r-1}) u^{|\lambda(\pm)_{z-1}|}}{(1+u) q^{[\sum_{h<i}
hm_h(\lambda(\pm)_{z-1})m_i(\lambda(\pm)_{z-1}) + \frac{1}{2} \sum_i
(i-1)m_i(\lambda(\pm)_{z-1})^2]} \prod_i A_O((z-1)^i)} \\
&&Prob(\Lambda(\pm)_{z+1}=\lambda(\pm)_{z+1})\\
 &=&  \frac{ \prod_{r=1}^{\infty}
(1-u^2/q^{2r-1}) u^{|\lambda(\pm)_{z+1}|}}{(1+u)q^{[\sum_{h<i}
hm_h(\lambda(\pm)_{z+1})m_i(\lambda(\pm)_{z+1}) + \frac{1}{2} \sum_i
(i-1)m_i(\lambda(\pm)_{z+1})^2]} \prod_i A_O((z+1)^i)} \\
Prob(\Lambda_{\phi}=\lambda_{\phi}) & = & \frac{ \prod_{r=1}^{\infty}
(1+(-1)^r \frac{u^{deg(\phi)}}{q^{deg(\phi)r/2}})u^{deg(\phi)
\cdot |\lambda_{\phi}|}}{q^{deg(\phi) [\sum_{h<i}
hm_h(\lambda_{\phi})m_i(\lambda_{\phi}) + \frac{1}{2} \sum_i
(i-1)m_i(\lambda_{\phi})^2]} \prod_i A_O(\phi^i)} \ if \
\phi=\bar{\phi} \neq z \pm 1\\
 Prob(\Lambda_{\phi}=\lambda_{\phi})
 & = & \frac{ \prod_{r=1}^{\infty}
(1- \frac{u^{2 deg(\phi)}}{q^{deg(\phi)r}}) u^{2 deg(\phi)
\cdot |\lambda_{\phi}|}}{q^{2 deg(\phi) [ \sum_{h<i}
hm_h(\lambda_{\phi})m_i(\lambda_{\phi}) + \frac{1}{2} \sum_i
(i-1)m_i(\lambda_{\phi})^2]} \prod_i A_O(\phi^i)} \ if \
\phi \neq \bar{\phi}.
\end{eqnarray*} Furthermore, setting $u=1$ in these formulas yields
the laws arising from the $n \rightarrow \infty$ limit of conjugacy
classes of a uniformly chosen element of $O(n,q)$, where the $+,-$
sign is chosen with probability $1/2$. The random variables
corresponding to different polynomials are independent, up to the fact
that $\Lambda_{\phi}= \Lambda_{\bar{\phi}}$. \end{theorem}

	For the same reasons as with the symplectic groups, the only
case remaining to be understood is the measure of the partition
corresponding to the polynomial $z-1$. Combining Theorem \ref{distSp}
with Lemma \ref{combinatorics} leads one to the following measure on
orthogonal signed partitions:

\[ M^{\pm}_{O,u}(\lambda(\pm)) = \frac{1}{1+u} \prod_{r=1}^{\infty}
(1-u^2/q^{2r-1}) \frac{u^{|\lambda(\pm)|}}{q^{1/2 [\sum_i
(\lambda(\pm)_i')^2 - \sum_i m_i(\lambda(\pm))^2 ]} \prod_i
A_{O}((z-1)^i)}. \] Forgetting about signs (i.e. lumping together
some conjugacy classes) yields a measure on underlying shapes which
will be denoted by $M_{O,u}$. Using the formulas for the sizes of the
finite symplectic and orthogonal groups given at the beginning of this
section, one arrives at the expression:

\begin{eqnarray*}
M_{O,u}(\lambda) & = & \frac{1}{1+u} \prod_{r=1}^{\infty} (1-u^2/q^{2r-1})
\frac{u^{|\lambda|}}{ q^{1/2 [\sum_i (\lambda_i')^2 - \sum_i
m_i^2 ]} \prod_{i=0 \ mod \ 2} (q^{\frac{m_i^2}{4}-\frac{m_i}{2}} \prod_{l=1}^{m_i/2} (q^{2l}-1))}\\
&  & \frac{1}{\prod_{i=1 \ mod \ 2 \atop m_i=0 \ mod \ 2} (q^{\frac{m_i^2}{4}-m_i} \prod_{l=1}^{m_i/2} (q^{2l}-1)) \prod_{i=0 \ mod \ 2 \atop m_i=1 \ mod \ 2} (q^{\frac{m_i^2-1}{4}} \prod_{l=1}^{(m_i-1)/2} (q^{2l}-1))}
\end{eqnarray*}

\section{Markov Chain Description} \label{Markov} 

	This section proves the Markov chain descriptions of conjugacy
classes as advertised in the introduction. The first goal is Theorem
\ref{Rudvalis}, which generalizes work of Rudvalis and Shinoda
\cite{RS} (proved by different methods and not in the language of
partitions). First, a lemma of Cauchy is needed.

\begin{lemma} \label{Cauchy} (\cite{A}, p.20) If $|q|>1$, \[
\prod_{m=0}^{\infty} (1-z/q^m)^{-1} = 1+\sum_{n=1}^{\infty}
\frac{z^n}{q^{n^2-n}(1-1/q)(1-1/q^2) \cdots(1-1/q^n) (1-z)(1-z/q)
\cdots (1-z/q^{n-1})}\] \end{lemma}

 Let $G$ denote either $Sp$ or $O$, and let $P_{G,u}(i)$ be the
probability that a partition chosen from the measure $M_{G,u}$ has $i$
parts. Let

\begin{eqnarray*}
P'_{Sp,u}(i) & = & \frac{P_{Sp,u}(i)}{\prod_{i=1}^{\infty} (1-u^2/q^{2i-1})}\\
P'_{O,u}(i) & = & \frac{(1+u) P_{O,u}(i)}{\prod_{i=1}^{\infty} (1-u^2/q^{2i-1})}.
\end{eqnarray*}

\begin{theorem} \label{Rudvalis}
 
\begin{eqnarray*}
P_{Sp,u}(2k) & = & \prod_{i=1}^{\infty} (1-u^2/q^{2i-1}) \frac{u^{2k}}{q^{2k^2+k}(1-u^2/q)(1-1/q^2) \cdots (1-u^2/q^{2k-1})(1-1/q^{2k})}\\
P_{Sp,u}(2k+1) & = & \prod_{i=1}^{\infty} (1-u^2/q^{2i-1}) \frac{u^{2k+2}}{q^{2k^2+3k+1}(1-u^2/q)(1-1/q^2) \cdots (1-1/q^{2k})(1-u^2/q^{2k+1})}\\
P_{O,u}(2k) & = & \frac{\prod_{i=1}^{\infty} (1-u^2/q^{2i-1})}{1+u} \frac{u^{2k}}{q^{2k^2-k}(1-u^2/q)(1-1/q^2) \cdots (1-u^2/q^{2k-1})(1-1/q^{2k})}\\
P_{O,u}(2k+1) & = & \frac{\prod_{i=1}^{\infty} (1-u^2/q^{2i-1})}{1+u} \frac{u^{2k+1}}{q^{2k^2+k}(1-u^2/q)(1-1/q^2) \cdots (1-1/q^{2k})(1-u^2/q^{2k+1})}.\\
\end{eqnarray*} \end{theorem}

\begin{proof} Using only the facts that $M_{Sp,u}$ and $M_{O,u}$
define a measure (i.e. not necessarily a probability measure), the
proofs of Theorems \ref{MarkovSp} and \ref{MarkovO} will establish the
equations:

\begin{eqnarray*}
P'_{Sp,u}(a) & = & \sum_{b \leq a \atop a-b \ even}
\frac{u^a P'_{O,u}(b)}{q^{(a^2-b^2+2(a+1)b)/4} (q^{a-b}-1) \cdots
(q^4-1)(q^2-1)}\\
P'_{O,u}(a) & = & \sum_{b \leq a \atop \ a-b \ even}
\frac{u^a P'_{Sp,u}(b) q^{(a-b)^2/4}}{q^{(a^2+b-2a)/2} (q^{a-b}-1) \cdots
(q^4-1)(q^2-1)}\\
&  & + \sum_{b \leq a \atop a-b \ odd} \frac{u^a P'_{Sp,u}(b) q^{((a-b)^2-1)/4}} {q^{(a^2-a)/2} (q^{a-b-1}-1) \cdots (q^4-1)(q^2-1)}.
\end{eqnarray*}

	To get a recurrence relation for the $P'_{Sp,u}(a)$'s, one
simply plugs the second equation into the first. Similarly one obtains
a recurrence relation for the $P'_{O,u}(a)$'s. These recurrences allow
one to solve for $P'_{G,u}(a)$ in terms of $P'_{G,u}(0)$, implying
that the formulas for $P_{G,u}(a)$ are proportional to the asserted
values. Thus it is enough to prove that the asserted formulas for
$P_{G,u}(a)$ satisfy the equation $\sum_{a \geq 0} P_{G,u}(a)=1$. This
follows readily from Lemma \ref{Cauchy}. \end{proof}

	Before continuing, we pause to indicate how the formulas of
Theorem \ref{Rudvalis} can be used to deduce group theoretic results
which are normally proved by techniques such as character theory and
Moebius inversion. The first set of results, Corollary
\ref{sympleccor}, considers only the symplectic groups. The same
technique would give results for the sum of $+,-$ type orthogonal
groups. Since the cycle index for the difference of orthogonal groups
also factors, one could rework all of the paper until now to give
measures corresponding to the difference of $+,-$ type orthogonal
groups, apply the same technique, and then average the results to get
theorems about groups over a given $+$ or $-$ type. This does not
deserve to be done publicly.

\begin{cor} \label{sympleccor}

\begin{enumerate}

\item (Steinberg, pg. 156 of \cite{H}) The number of unipotent
elements in $Sp(2n,q)$ is $q^{2n^2}$.

\item (\cite{RS}) The probability that an randomly chosen element of
$Sp(2n,q)$ has a $2k$ dimensional fixed space is

\[ \frac{1}{|Sp(2k,q)|} \sum_{i=0}^{n-k} \frac{(-1)^i (q^2)^{i \choose
2}}{|Sp(2i,q)| q^{2ik}}.\] The probability that an randomly chosen
element of $Sp(2n,q)$ has a $2k+1$ dimensional fixed space is \[
\frac{1}{|Sp(2k,q)| q^{2k+1}} \sum_{i=0}^{n-k-1} \frac{(-1)^i (q^2)^{i
\choose 2}}{|Sp(2i,q)| q^{2i(k+1)}}.\]

\end{enumerate}

\end{cor}

\begin{proof} The arguments are completely analogous to those for
$GL(n,q)$ in \cite{F2} (Corollary 1 and Theorem 6), using the cycle
index of the finite symplectic groups \cite{F1}. \end{proof}

	Rudvalis and Shinoda (loc. cit.) also considered the
probability that the fixed space of a random element of a finite
classical group has a given isometry type. For the finite unitary
groups, the isometry classes are parameterized by pairs $(s,t)$ of
natural numbers such that $s+2t \leq n$. Here a subspace $W$ of $V$
has type $(s,t)$ if $dim(W/rad(W))=s$ and $dim(rad(W))=t$. Theorem
\ref{unitarycor} uses cycle index techniques to give new proofs of
their results for the finite unitary groups. Exactly the same methods
work for the finite symplectic and orthogonal groups, but we spare the
reader the details.

\begin{cor} \label{unitarycor} The probability that an element of
$U(n,q)$ has isometry type corresponding to the pair $(s,t)$ is

\[ \frac{ \sum_{i=0}^{n-2s-t} \frac{(-1/q)^{(t+1)i} (-1/q)^{{i \choose
2}}}{(1+1/q)(1-1/q^2) \cdots (1-(-1)^i/q^i)}}{q^{s^2+2st}
(1+1/q)(1-1/q^2) \cdots (1-(-1)^s/q^s) (1+1/q)(1-1/q^2) \cdots
(1-(-1)^t/q^t)}.\] In the $n \rightarrow \infty$ limit, this converges
to \[ \frac{ \prod_{r=0}^{\infty}
(\frac{1}{1+1/q^{2r+1}})}{q^{s^2+2st} (1+1/q)(1-1/q^2) \cdots
(1-(-1)^s/q^s) (1-1/q^2)(1-1/q^4) \cdots (1-1/q^{2t})}.\] \end{cor}

\begin{proof} The most important observation (see \cite{FNP} for a
readable proof) is that the fixed space of an element $\alpha$ of
$U(n,q)$ has isometry type $(s,t)$ precisely when the partition
corresponding to the polynomial $z-1$ in the rational canonical form
of $\alpha$ satisfies $\lambda_1'=s+t$ and $\lambda_2'=t$. In other
words, the partition has $s+t$ parts and $s$ $1$'s. Let $[u^n] f(u)$
denote the coefficient of $u^n$ in some polynomial $f(u)$. Then one uses
the cycle index for the unitary groups as in \cite{F2} to see
that the sought probability for $U(n,q)$ is

\[ [u^n] \frac{1}{1-u} \sum_{\lambda: \lambda_1'=s+t, \lambda_2'=t}
M_{U,u}(\lambda).\] Using the fact that 

\[ M_{U,u}(\lambda) = \prod_{r=1}^{\infty} (1+u/(-q)^r)
\frac{u^{|\lambda|}}{q^{\sum_i (\lambda_i')^2} \prod_i
(1+1/q)(1-1/q^2) \cdots (1+(-1)^{m_i+1}/q^{m_i})},\] this becomes

\begin{eqnarray*}
&  & [u^n] \frac{u^{s+t}}{(1-u)q^{(s+t)^2} (1+1/q)(1-1/q^2) \cdots
(1-(-1)^s/q^s)} \sum_{\lambda: \lambda_1'=t} M_{U,u}(\lambda)\\
& = & \frac{1}{q^{(s+t)^2} (1+1/q)(1-1/q^2) \cdots
(1-(-1)^s/q^s)} [u^{n-s-t}]  \frac{1}{1-u} \sum_{\lambda: \lambda_1'=t} M_{U,u}(\lambda)\\
& = &  \frac{\sum_{i=0}^{n-2s-t}
\frac{(-1/q)^{(t+1)i} (-1/q)^{{i \choose 2}}}{(1+1/q)(1-1/q^2) \cdots
(1-(-1)^i/q^i)}}{q^{s^2+2st} (1+1/q)(1-1/q^2) \cdots (1-(-1)^s/q^s)
(1+1/q)(1-1/q^2) \cdots (1-(-1)^t/q^t)},
\end{eqnarray*} where the last equality is in the proof of Theorem 6 in \cite{F2}.

	The formula for the $n \rightarrow \infty$ limit follows from
the well-known identity

\[ \prod_{r=1}^{\infty} (1-v/w^r) = \sum_{n=0}^{\infty}
\frac{(-v)^n}{(w^n-1) \cdots (w-1)}.\] Alternatively, it follows from
the principle that the limit as $n \rightarrow \infty$ of $f(u)/(1-u)$
is $f(1)$ if $f$ has a Taylor expansion around $0$ converging in a
circle of radius $1$, together with the formula for 

\[ \sum_{\lambda: \lambda_1'=t} M_{U,1}(\lambda) \] given in Theorem 5
of \cite{F2}. \end{proof}

	Lemma \ref{recur} is crucial and motivated the combinatorial
moves made in rewriting the formulas for $M_{Sp,u}$ and $M_{O,u}$ in
Section \ref{recall}. In all that follows $M_{G,x}$ will denote the
probability of an event $X$ under the measure $M_{G,u}$ with $G$ equal
to $Sp$ or $O$.

\begin{lemma} \label{recur}
\begin{enumerate}
\item If $i$ is odd then

\[ M_{Sp,u}(\lambda_1'=s_1,\cdots,\lambda_{i-1}'=s_{i-1},\lambda_i'=k) =
 \frac{\prod_{r=1}^{\infty} (1-u^2/q^{2r-1}) u^{s_1+
\cdots+s_{i-1}}  P_{Sp,u}'(k)} {q^{(s_1^2+\cdots
+s_{i-1}^2-m_1^2-\cdots-m_{i-1}^2)/2} \prod_{j=1}^{i-1}
A_{Sp}((z-1)^{m_j})}. \]

\item If $i$ is even then

\begin{eqnarray*}
&&M_{Sp,u}(\lambda_1'=s_1,\cdots,\lambda_{i-1}'=s_{i-1},\lambda_i'=k)\\
& = &
 \frac{\prod_{r=1}^{\infty} (1-u^2/q^{2r-1}) u^{s_1+
\cdots+s_{i-1}} P_{O,u}'(k)} {q^{(k+s_1^2+\cdots
+s_{i-1}^2-m_1^2-\cdots-m_{i-1}^2)/2} \prod_{j=1}^{i-1}
A_{Sp}((z-1)^{m_j})}.
\end{eqnarray*}

\item If $i$ is odd then

\begin{eqnarray*}
&& M_{O,u}(\lambda_1'=s_1,\cdots,\lambda_{i-1}'=s_{i-1},\lambda_i'=k)\\
& = &
 \frac{ \prod_{r=1}^{\infty} (1-u^2/q^{2r-1})u^{s_1+
\cdots+s_{i-1}} P_{O,u}'(k)} {(1+u) q^{(s_1^2+\cdots
+s_{i-1}^2-m_1^2-\cdots-m_{i-1}^2)/2} \prod_{j=1}^{i-1}
A_{O}((z-1)^{m_j})}.
\end{eqnarray*}

\item If $i$ is even then

\begin{eqnarray*}
&& M_{O,u}(\lambda_1'=s_1,\cdots,\lambda_{i-1}'=s_{i-1},\lambda_i'=k)\\
& = &
 \frac{ \prod_{r=1}^{\infty} (1-u^2/q^{2r-1}) u^{s_1+
\cdots+s_{i-1}} q^{k/2}  P_{Sp,u}'(k) } {(1+u) q^{(s_1^2+\cdots
+s_{i-1}^2-m_1^2-\cdots-m_{i-1}^2)/2} \prod_{j=1}^{i-1}
A_{O}((z-1)^{m_j})}.
\end{eqnarray*}

\end{enumerate}
\end{lemma}

\begin{proof} The idea for all of the proofs is the same; hence we
prove part two as follows:

\begin{eqnarray*}
& & M_{Sp,u}(\lambda_1'=s_1,\cdots,\lambda_{i-1}'=s_{i-1},\lambda_i'=k)\\
& = & \frac{ \prod_{r=1}^{\infty} (1-u^2/q^{2r-1}) u^{s_1+
\cdots+s_{i-1}}} {q^{(s_1^2+\cdots
+s_{i-1}^2-m_1^2-\cdots-m_{i-1}^2)/2} \prod_{j=1}^{i-1}
A_{Sp}((z-1)^{m_j})}\\
&& \cdot \sum_{\lambda_i'=k \geq \lambda_{i+1}' \geq \cdots \geq 0} \frac{u^{\sum_{j \geq i} \lambda_j'}}{q^{\sum_{j
\geq i} (\lambda_j'^2-m_j^2)/2} \prod_{j \geq i} A_{Sp}((z-1)^{m_j})}\\
& = & \frac{ \prod_{r=1}^{\infty} (1-u^2/q^{2r-1}) u^{s_1+
\cdots+s_{i-1}}} {q^{(k+s_1^2+\cdots
+s_{i-1}^2-m_1^2-\cdots-m_{i-1}^2)/2} \prod_{j=1}^{i-1}
A_{Sp}((z-1)^{m_j})}\\
&& \cdot \sum_{\lambda_i'=k \geq \lambda_{i+1}' \geq \cdots \geq 0} \frac{u^{\sum_{j \geq i} \lambda_j'}}{q^{\sum_{j
\geq i} (\lambda_j'^2-m_j^2)/2} \prod_{j \geq i} A_{O}((z-1)^{m_j})}\\
& = & \frac{ \prod_{r=1}^{\infty} (1-u^2/q^{2r-1}) u^{s_1+
\cdots+s_{i-1}}} {q^{(k+s_1^2+\cdots
+s_{i-1}^2-m_1^2-\cdots-m_{i-1}^2)/2} \prod_{j=1}^{i-1}
A_{Sp}((z-1)^{m_j})}\\
&& \cdot \sum_{\lambda_1'=k \geq \lambda_{2}' \geq \cdots \geq 0} \frac{u^{\sum_{j \geq 1} \lambda_j'}}{q^{\sum_{j
\geq 1} (\lambda_j'^2-m_j^2)/2} \prod_{j \geq 1} A_{O}((z-1)^{m_j})}\\
& = & \frac{ \prod_{r=1}^{\infty} (1-u^2/q^{2r-1}) u^{s_1+
\cdots+s_{i-1}}} {q^{(k+s_1^2+\cdots
+s_{i-1}^2-m_1^2-\cdots-m_{i-1}^2)/2} \prod_{j=1}^{i-1}
A_{Sp}((z-1)^{m_j})} P_{O,u}'(k),
\end{eqnarray*} as desired. Note that the meat of the lemma is the second equality, which follows from the formulas for $A_{Sp}$ and $A_O$. The third equality is simply a relabelling of subscripts.
\end{proof}

	For the statements of Theorems \ref{MarkovSp} and
\ref{MarkovO}, we define two Markov chains on the integers. Recall
that
 
\begin{eqnarray*}
P'_{Sp,u}(2k) & = & \frac{u^{2k}}{q^{2k^2+k}(1-u^2/q)(1-1/q^2) \cdots (1-u^2/q^{2k-1})(1-1/q^{2k})}\\
P'_{Sp,u}(2k+1) & = & \frac{u^{2k+2}}{q^{2k^2+3k+1}(1-u^2/q)(1-1/q^2) \cdots (1-1/q^{2k})(1-u^2/q^{2k+1})}\\
P'_{O,u}(2k) & = & \frac{u^{2k}}{q^{2k^2-k}(1-u^2/q)(1-1/q^2) \cdots (1-u^2/q^{2k-1})(1-1/q^{2k})}\\
P'_{O,u}(2k+1) & = & \frac{u^{2k+1}}{q^{2k^2+k}(1-u^2/q)(1-1/q^2) \cdots (1-1/q^{2k})(1-u^2/q^{2k+1})}.\\
\end{eqnarray*} The chains $K_1,K_2$ are defined on the natural numbers with
transition probabilities

\[ K_1(a,b) =             \left\{ \begin{array}{ll}
	\frac{u^aP'_{O,u}(b)}{P'_{Sp,u}(a) q^{\frac{a^2-b^2+2(a+1)b}{4}} (q^{a-b}-1) \cdots (q^4-1)(q^2-1)}             & \mbox{if $a-b$ even}\\
												0	& \mbox{if
$a-b$ odd}
																																				\end{array}
			\right.			 \]

\[ K_2(a,b) =  \left\{ \begin{array}{ll}
	\frac{u^aP'_{Sp,u}(b) q^{(a-b)^2/4}}{P'_{O,u}(a) q^{\frac{a^2+b}{2}-a} (q^{a-b}-1) \cdots (q^4-1)(q^2-1)}             & \mbox{if $a-b$ even}\\
											\frac{u^aP'_{Sp,u}(b) q^{((a-b)^2-1)/4}}{P'_{O,u}(a) q^{\frac{a^2-a}{2}} (q^{a-b-1}-1) \cdots (q^4-1)(q^2-1)}    & \mbox{if
$a-b$ odd}
													\end{array}
			\right.  \] The fact that these transition
probabilities add up to one will follow from the proof of Theorem
\ref{MarkovSp}.

\begin{theorem} \label{MarkovSp} Starting with $\lambda_1'$
distributed as $P_{Sp,u}$, define in succession
$\lambda_2',\lambda_3',\cdots$ according to the rules that if
$\lambda_i'=a$, then $\lambda_{i+1}'=b$ with probability $K_1(a,b)$ if
$i$ is odd and probability $K_2(a,b)$ if $i$ is even. The resulting
partition is distributed according to $M_{Sp,u}$. \end{theorem}

\begin{proof} The $M_{Sp,u}$ probability of choosing a partition
with $\lambda_i'=s_i$ for all $i$ is

\[ M_{Sp,u}(\lambda_1'=s_1) \prod_{i=1}^{\infty}
\frac{M_{Sp,u}(\lambda_1'=s_1,\cdots,\lambda_{i+1}'=s_{i+1})}
{M_{Sp,u}(\lambda_1'= s_1,\cdots,\lambda_{i}'=s_{i})}.\] Since
$M_{Sp,u}(\lambda_1'=s_1)$ is equal to $P_{Sp,u}(s_1)$ by definition,
it it is enough to prove two claims: first that for every choice of
$i,a,b,s_1,\cdots,s_{i-1}$,

\[ \frac{M_{Sp,u}(\lambda_1'=s_1,\cdots,
\lambda_{i-1}'=s_{i-1},\lambda_i'=a,\lambda_{i+1}'=b)}{M_{Sp,u}(
\lambda_1'=s_1,\cdots, \lambda_{i-1}'=s_{i-1},\lambda_i'=a)} \] is
equal to the asserted transition rule probability for moving from
$\lambda_i'=a$ to $\lambda_{i+1}'=b$, and second that the transition rule
probabilities sum to one.

	The first claim follows from Lemma \ref{recur}. For the second
claim, observe that \[ \sum_{b \leq a}
\frac{M_{Sp,u}(\lambda_1'=s_1,\cdots,
\lambda_{i-1}'=s_{i-1},\lambda_i'=a,\lambda_{i+1}'=b)}{M_{Sp,u}(
\lambda_1'=s_1,\cdots, \lambda_{i-1}'=s_{i-1},\lambda_i'=a)}
 =1,\] because $M_{Sp,u}$ is a measure and the columns of a partition
are non-increasing in size as one moves to the right. Since $\sum_{i
\geq 0} P_{Sp,u}(i)=1$, it follows that $M_{Sp,u}$ is a probability
measure, as promised earlier.  \end{proof}

	Theorem \ref{MarkovO} gives the analogous result for the
orthogonal groups. As the proof method is the same as for the
symplectic groups, we merely record the result.

\begin{theorem} \label{MarkovO} Starting with $\lambda_1'$ distributed
as $P_{O,u}$, define in succession $\lambda_2',\lambda_3',\cdots$
according to the rules that if $\lambda_i'=a$, then $\lambda_{i+1}'=b$
with probability $K_2(a,b)$ if $i$ is odd and $K_1(a,b)$ if $i$ is
even. The resulting partition is distributed according to
$M_{O,u}$. \end{theorem}

	We close the paper with the following remarks.

{\bf Remarks:}
\begin{enumerate}

\item Theorem \ref{MarkovSp} and Theorem \ref{MarkovO} allow one to
draw exact samples from the measures $M_{Sp,u}$ or $M_{O,u}$. First
recall that sampling from discrete distributions $P$ with known
formulas is straightforward; simply pick $U$ uniformly in $[0,1]$ and
find the value of $j$ such that $ \sum_{i=0}^j P(i) < U <
\sum_{i=0}^{j+1} P(i)$. This allows one to sample from $P_{Sp,u}$ or
$P_{O,u}$. Then move according to the appropriate Markov chains.

\item For the case of $M_{Sp,u}$, one can view the algorithm of
Theorem \ref{MarkovSp} slightly differently. One starts with an
imaginary $0$th column of size approaching infinity, and then gets
$\lambda_1'$ by transitioning according to the chain $K_2$. It is
straightforward to verify that the resulting distribution of the first
column size agrees with $P_{Sp,u}$. This viewpoint was useful in the
general linear and unitary cases \cite{F4}.

\item As noted in the introduction, it is possible that the measures
$M_{Sp,u}$ and $M_{O,u}$ are related to generalizations of the
Rogers-Ramanujan identities, in analogy with the corresponding
measures for $GL(n,q)$. In this regard observe that

\begin{eqnarray*}
 \sum_{\lambda:\lambda_2'=0} M_{Sp,u}(\lambda) &= &
\prod_{r=1}^{\infty} (1-u^2/q^{2r-1}) \sum_{n=0}^{\infty}
\frac{u^{2n}}{q^{n^2}(q^{2n}-1)\cdots(q^2-1)}\\
 & = & \prod_{r=1}^{\infty} (1-u^2/q^{2r-1}) \sum_{n=0}^{\infty}
\frac{u^{2n}}{q^{2n^2+n}(1-1/q^{2})\cdots(1-1/q^{2n})}.
\end{eqnarray*} For the value $u=q^{-1/2}$, the sum has a product
expansion by a Rogers Ramanujan identity.

\item Although a probabilistic understanding of $M_{Sp,u}$ and
$M_{O,u}$ has been given, a viewpoint explaining the products in
Theorem \ref{Rudvalis} in terms of certain random variables being
independent would be desirable. This was possible for the finite
general linear and unitary groups \cite{F2}.

\item As mentioned in the introduction, the clumping of conjugacy
classes given by looking at the underlying shape was not necessary. To
modify things to take signs into account, $K_1$ starts at a number $a$
but outputs an ordered pair $(b,\pm)$, that is a choice of sign
associated with $b$. In the $+$ (resp. $-$) case, the expression
$P_{O,u}'(b)$ in the numerator of the definition of $K_1$ is replaced
by the probability that the partition under the measure $M_{O,u}$ has
$b$ parts and a $+$ (resp. $-$) choice for the parts of size
1. Unfortunately we do not know of simple formulas for these
probabilities analogous to Theorem \ref{Rudvalis}. Formulas can be
inferred from the paper \cite{RS}, but the result involves unpleasant
sums and does not seem useful. The transition probabilities for $K_2$
are also affected: the $P_{O,u}'(a)$ are replaced the same way as for
$K_1$, and letting $\epsilon$ be the sign associated to $a$, the
transition probabilities are multiplied by an additional factor of \[
\frac{1}{|O^{\epsilon}(a-b,q)|}
\frac{1}{\frac{1}{|O^+(a-b,q)|}+\frac{1}{O^-(a-b,q)|}}.\] It is not
necessarily surprising that the theory is nicer when conjugacy classes
are lumped; the cycle index only factors for sums or differences of
orthogonal groups.

\end{enumerate}

\section*{Acknowledgements} This work is a natural follow-up to the
author's Ph.D. thesis; the author expresses deep gratitude to Persi
Diaconis for having introduced him to this lovely part of
mathematics. A summer of conversations with Peter M. Neumann and
Cheryl E.  Praeger broadened the author's understanding of conjugacy
classes. This work received the financial support of an NSF
Postdoctoral Fellowship.

\end{document}